\documentclass{commat}

\title{Summability characterizations of positive sequences}

\author{Douglas Azevedo and Thiago P. de Andrade}

\affiliation{
	\address{Universidade Tecnologica Federal do Parana.  Av. Alberto Carazzai, n 1640. CEP 863000-000, Caixa Postal 238, Cornelio Procopio, Parana, Brazil}
	\email{douglasa@utfpr.edu.br, thiagopinguello@hotmail.com}
	}

\abstract{
In this paper, we propose extensions for the classical Kummer's test, which is a very far-reaching criterion that provides sufficient and necessary conditions for convergence and divergence of series of positive terms. Furthermore, we present and discuss some interesting consequences and examples such as extensions of Olivier's theorem and Raabe, Bertrand and Gauss's test.	
}



\firstpage{81}

\VOLUME{30}

\DOI{https://doi.org/10.46298/cm.9290}

\begin{paper}

\section{Introduction}

The Kummer's test is an advanced theoretical test which provides necessary and sufficient conditions that ensures convergence and divergence of series of positive terms. Below we present the statement of this result. Its proof and some additional historical background may be found in~\cite{Ludmila}, \cite{Knopp}, \cite{Tong:2004}.

\begin{theorem}
(Kummer's test)\label{Kummer}
Consider the series \(\sum a_{n}\) where \(\{a_{n}\}\) is a~sequence of positive real numbers.
\begin{enumerate}
\item[(i)] The series \(\sum a_{n}\) converges if and only if there exist a~sequence \(\{q_{n}\}\), a~real number \(c>0\) and an integer \(N\geq 1\) for which
\[
q_{n}\frac{a_{n}}{a_{n+1}}-q_{n+1}\geq c, \qquad n \geq N.
\]

\item[(ii)] The series \(\sum a_{n}\) diverges if, and only if there exist a~sequence \(\{q_{n}\}\) and an integer \(N\geq 1\) for which
\(\sum \frac{1}{q_{n}}\) is a~divergent series and
\[
q_{n}\frac{a_{n}}{a_{n+1}}-q_{n+1}\leq 0, \qquad n \geq N.
\]
\end{enumerate}
\end{theorem}

 Besides providing an extremely far-reaching characterization of convergence and divergence of series with positive terms, the importance of Kummer's test it is mostly ratified by its implications. For instance, Bertrand's test, Gauss's test, Raabe's test~\cite{Tong:2004} are all special cases of Theorem~\ref{Kummer}. Kummer's test may be also usefull to characterize convergence in normed vector spaces~\cite[p. 7]{Muscat} and applications of this test
 can be found in other branches of Analysis, such as difference equations~\cite{Gyori}, as well.

 On the other hand, turning our focus to series of the form \(\sum c_n a_n\), there are only few results dealing with this type of series. The Abel's test and test of Dedekind and Du-Bois Reymond (see for instance,~\cite[p. 315]{Knopp},~\cite{Hadamard}) are probably the most famous, since they deal with general series of complex numbers. These tests provide conditions that ensure convergence by means of independent assumptions on \(\{c_{n}\}\) and \(\{a_{n}\}\). In this context, the main feature of our results (Theorem~\ref{thm1} and Theorem~\ref{thm2}) is that they characterize the relation between the sequences \(\{a_{n}\}\) and \(\{c_{n}\}\) in order to ensure necessary and sufficient conditions for the convergence and divergence of the series \(\sum c_n a_n\), respectively. Moreover, we present some examples and interesting consequences of this characterization. In particular, generalized versions of Raabe's, Bertrand's and Gauss's test for convergence and divergence of series of the form \(\sum c_n a_n\)
are obtained. Another important consequence of Theorem~\ref{thm1} is that it is possible to
show that Olivier's theorem (see, for instance~\cite[p. 124]{Knopp} or~\cite{Const}, \cite{Salat} for more information)
still holds when the monotonicity assumption on the sequence of positive terms \(\{a_n\}\)
is replaced by an additional assumption on a~auxiliary sequence. We also present consequences of Theorem~\ref{thm1} when it is combined to the well-know Abel summation formula
 and the Cauchy condensation theorem. We refer to~\cite[pp. 120 and 313]{Knopp} for more details on these results.

 The rest of the paper is organized as follows.
 In Section~\ref{subsum}, we present necessary and sufficient conditions for convergence/divergence of series generated by subsequences by extending Theorem~\ref{Kummer}.
 In Section~\ref{cnan} we present the results dealing with convergence and divergence of series of the form \(\sum c_n a_n\). The main ideia is to obtain necessary and sufficient conditions by means of an extension of Theorems~\ref{conv} and Theorem~\ref{div}. As we show, we characterize the relation between the sequences \(\{c_n\}\) and \(\{a_n\}\) that ensures convergence and divergence of the series.
 In Section~\ref{EC} we present some consequences of the results obtained.

\section{An extension of Kummer's test: I}\label{subsum}

In this section we present a~first extension of Theorem~\ref{Kummer}. Its main feature is that it showns that is possible to obtain information about the summability of a~sequence of positive real numbers based on the relation between non-consecutive elements of this sequence. In partiular, the idea is to characterize the summability of a~sequence by comparing it to the elements of the translated sequence \(\{a_{n+m}, \ n\geq 1\}\), for some~\(m \ge 1\).

The first main result of this section is presented below.

\begin{theorem}\label{conv}
Let \(\{a_n\}\) be a~sequence of positive real numbers and \(m\geq 1\) any
fixed positive integer.

If there exists a~positive sequence \(\{q_n\}\) such that
\[
q_n\frac{a_n}{a_{n+m}}-q_{n+m}\geq c,
\]
for some \(c>0\), for all \(n\) sufficiently large, then \(\sum a_{n}\) converges.
The converse holds as well.
\end{theorem}
\begin{proof}
 From the assumption we get that
\[
q_n a_n-a_{n+m}q_{n+m}\geq ca_{n+m},
\]
for all \(n>N\), for some \(N\) large.
Hence
\[
\sum_{n = N+1}^{N+k} q_n a_n-a_{n+m}q_{n+m}\geq c \sum_{n = N+1}^{N+k} a_{n+m},
\]
for all \(k\geq 1\). That is, by the telescopic sum and considering without loss of generality \(k > m\), we have
\begin{align*}
q_{N+1}a_{N+1} + \dots + q_{N+m}a_{N+m} - a_{N+k+1}q_{N+k+1} - \dots -{ }&{ }a_{N+k+m}q_{N+k+m} \\
&{ }\geq c \sum_{n = N+1}^{N+k} a_{n+m},
\end{align*}
for all \(k>m\). Since \(\{a_n\}\) and \(\{q_n\}\) are positive, the left side of previous inequality is less than \(q_{N+1} a_{N+1}+ \dots + q_{N+m}a_{N+m}\) and then the series
\(\sum a_{n+m}\) converges. Therefore, \(\sum a_{n}\) also converges.

Conversely, if \(\sum a_n\) converges, \(\sum a_n = S\) say,
then let us write \(\sum a_{n+m-1} = S_m\), for \(m\geq 1\), positive integer.
Let us define \(\{q_n\}\) as
\[
q_n = \frac{S_{m}-\sum_{i = 1}^{n}a_{i+m-1}}{a_{n}}, \qquad n = 1, 2, 3, \dots ,
\]
 thus, for this \(\{q_n\}\) we have that
\begin{align*}
q_n \frac{a_n}{a_{n+m}}-q_{n+m}
{ }&{ }= \frac{\sum_{i = n+1}^{n+m}a_{i+m-1}}{a_{n+m}}\\
{ }&{ }= 1+\frac{a_{n+m+1}+ \dots +a_{n+2m-1}}{a_{n+m}}\\
{ }&{ }> 1,
\end{align*}
for all \(n\geq 1\).
The proof is concluded.
\end{proof}

We proceed by presenting a~divergence version for the previous theorem.

\begin{theorem}\label{div}
Let \(\{a_n\}\) be a~sequence of positive real and \(m\geq 1\) a~fixed positive integer. If there exists a~positive sequence \(\{q_n\}\) such that \(\sum \frac{1}{q_{n}}\) diverges, \(q_n a_n\geq c>0\),
 and
\[
q_n\frac{a_n}{a_{n+m}}-q_{n+m}\leq 0,
\]
 for all \(n\) sufficiently large,
then
\(\sum_{n = 1}^{\infty} a_{n}\) diverges.
The converse holds, as well.
\end{theorem}
\begin{proof}
From the assumptions we obtain that there exists \(N>0\) such that
\[
q_n\frac{a_n}{a_{n+m}}-q_{n+m}\leq 0,
\]
 for all \(n\geq N\).
 As so,
\[
c\frac{1}{q_{n+m}}\leq a_{n+m},
\]
 for all \(n>N\). Since \(\sum 1/q_n\) diverges, we obtain from the comparsion test that \(\sum a_n\) diverges.

Conversely, suppose that \(\sum a_{n}\) diverges. Define for each \(n\geq 1\)
\[
q_n = \frac{\sum_{i = 1}^{n}a_i}{a_n}.
\]
Note that the definition implies \(q_1 = 1\), hence \(a_n q_n = \sum_{i = 1}^{n}a_{i}\geq a_1\), for all \(n\geq 1\), that is, \(a_n q_n\geq a_1 q_1>0\) for all \(n\geq 1\).
Clearly
\[
q_{n}\frac{a_{n}}{a_{n+m}}-q_{n+m}\leq 0,
\]
for all \(n\geq 1\).

Let us now show that \(\sum \frac{1}{q_{n}}\)
diverges. From the divergence of \(\sum a_{n}\), given any positive integer \(k\) there exists a~positive integer \(n\geq k\) such that
\begin{equation}\label{aux11}
 a_{k}+\dots+a_{n}\geq a_1+\dots+a_{k-1}.
\end{equation}
Due to~\eqref{aux11},
\begin{align*}
\sum_{j = k}^{n}\frac{1}{q_{j}}
{ }&{ }= \frac{a_{k}}{a_{1}+ \dots +a_{k}}+\dots+\frac{a_{n}}{a_{1}+ \dots +a_{n}}\\
{ }&{ }\geq \frac{a_{k}}{a_{1}+ \dots +a_{n}}+\dots+\frac{a_{n}}{a_{1}+ \dots +a_{n}}\\
{ }&{ }= \frac{1}{\frac{a_{1}+ \dots +a_{k-1}}{a_{k}+ \dots +a_{n}} +1}\\
{ }&{ }> \frac{1}{2}.
\end{align*}
Hence, \(\sum_{j = 1}^{n}\frac{1}{q_{j}}\) is not a~Cauchy sequence. Therefore the series \(\sum \frac{1}{q_{n}}\) diverges.
\end{proof}

\section{Extension of Kummer's test: II}\label{cnan}

Let us now turn our atention to series of the form \(\sum c_{n}a_{n}\) with positive terms.
The central idea in the following result is that it characterizes
the relation between the sequences \(\{c_{n}\}\) and \(\{a_{n}\}\) in order to ensure the convergence of the series.
The reader will note that the proof follows the same lines as the proof of Theorem~\ref{conv} and also, that it could be obtained
by some changes in the proof of Theorem~\ref{Kummer}, nevertherless, as the reader will also note, our proof
provides important informations about the relation between the sequences \(\{a_n\}\) and \(\{c_n\}\).

\begin{theorem}\label{thm1}
Consider the series \(\sum c_{n}a_{n}\) with \(\{a_{n}\}\) \(\{c_{n}\}\) sequences of positive real numbers.
 The series \(\sum c_{n}a_{n}\) converges if and only if that there exist a~sequence \(\{q_{n}\}\) of positive real numbers and a~positive integer \(N\geq 1\) for which
\[
q_{n}\frac{a_{n}}{a_{n+1}}-q_{n+1}\geq c_{n+1}, \quad n\geq N.
\]
\end{theorem}

\begin{proof}
 Let us show that \(\sum c_{n}a_{n}\) converges. For this, note that the condition
\[
q_{n}\frac{a_{n}}{a_{n+1}}-q_{n+1}\geq c_{n+1}, \quad n\geq N
\]
implies that
\begin{equation}\label{auxx}
a_{n}q_{n}\geq a_{n+1}(q_{n+1}+ c_{n+1}), \quad n\geq N.
\end{equation}
That is,
\begin{align*}
a_{N}q_{N}
{ }&{ }\geq a_{N+1}(q_{N+1}+c_{N+1})\\
{ }&{ }\geq a_{N+2}(q_{N+2}+c_{N+2})+a_{N+1}c_{N+1}\\
{ }&{\ \,}\vdots\\
{ }&{ }\geq a_{N+k}q_{N+k} + \sum_{i = 1}^{k}c_{N+i}a_{N+i}\\
{ }&{ }\geq \sum_{i = 1}^{k}c_{N+i}a_{N+i}>0,
\end{align*}
for all integer \(k\geq 0\). This implies the convergence of \(\sum c_{n}a_{n}\).

 For the converse, suppose that
\(S:= \sum c_{n}a_{n}\) and
let us define
\begin{equation}\label{pn}
q_{n} = \,\frac{S-\sum_{i = 1}^{n}c_{i}a_{i}}{a_{n}}, \quad n\geq N.
\end{equation}
 For this \(\{q_{n}\}\), clearly \(q_{n}>0\) for all \(n\geq 1\) and it is easy to check that
\[
q_{n}\frac{a_{n}}{a_{n+1}}-q_{n+1} = c_{n+1}, \quad n\geq N.
\qedhere
\]
\end{proof}

Some remarks:
\begin{enumerate}
\item[(i)] One can observe that it is, of course, possible to
reduce any series to this form, as any number can be
expressed as the product of two other numbers. Success in applying the above
theorem will depend on the skill with which the terms are so split up.

\item[(ii)] Note that in the first part of Theorem~\ref{thm1}, the assumption of
positivity of the sequences \(\{a_{n}\}\) and \(\{c_{n}\}\)
can be replaced by the following assumptions: \(\{a_{n}\}\) is positive and \(\{c_n\}\) is such that \(\sum_{i = 1}^{k}c_{i}a_{i}>0 \) for all
\(k\) sufficiently large.
\end{enumerate}

Next, we presente a~version of Kummer's test for divergent series
of the form \(\sum c_{n}a_{n}\).
The reader will note that it is more restrictive when
it is compared to Theorem~\ref{Kummer}-\((ii)\) however it may be suitable in some cases.

\begin{theorem}\label{thm2}
Consider the series \(\sum c_{n}a_{n}\) with \(\{a_{n}\}\) \(\{c_{n}\}\) sequences of positive real numbers.
\begin{enumerate}
\item[(i)] Suppose that
 there exist a~sequence \(\{q_{n}\}\) and a~positive integer \(N\) for which
\[
q_{n}\frac{a_{n}}{a_{n+1}}-q_{n+1}\leq -c_{n+1}, \quad n\geq N
\]
 with \(\sum \frac{1}{q_{n}}\) being a~divergent series. Then \(\sum a_{n}\), \(\sum \frac{1}{c_{n}}\), \(\sum (q_{n}-c_{n})a_{n}\) and \(\sum q_{n}a_{n}\) diverge. If, in addition, \(\sum \frac{c_{n}}{q_{n}}\) diverges
 then \(\sum c_{n}a_{n}\) diverges.

\item[(ii)] Suppose that both series \(\sum c_{n}a_{n}\) and \(\sum a_{n}\) diverge. Also,
 suppose that for every \(m\in\mathbb{N}\) there exists \(r\geq m\), \(r\in\mathbb{N}\), such that
\[
a_{m}+\dots+a_{r}\geq c_{m}a_{m}+\dots+c_{r}a_{r}.
\]
 Then there exist a~sequence \(\{q_{n}\}\) and a~positive integer \(N\geq 1\) such that
\(\sum\frac{1}{q_{n}}\) diverges
and
\[
q_{n}\frac{a_{n}}{a_{n+1}}-q_{n+1}\leq -c_{n+1}, \quad n\geq N.
\]
\end{enumerate}
\end{theorem}
\begin{proof}
To prove \((i)\) note that
\(\{q_{n}\}\) satisfies
\begin{gather}
a_{n+1}\geq \frac{q_{n}a_{n}}{q_{n+1}-c_{n+1}}, \quad n\geq N,\label{aux111} \\
0<q_{n+1}-c_{n+1}<q_{n+1},\quad n\geq N \label{aux1} \\
\textup{and} \nonumber \\
0<c_{n+1}<q_{n+1}, \quad n\geq N. \nonumber
\end{gather}

By last inequality and comparsion test we see that \(\sum \frac{1}{c_{n}}\) diverges. Next, using~\eqref{aux111} successively we see that
\begin{gather*}
a_{N+1}\geq \frac{q_{N}a_{N}}{q_{N+1}-c_{N+1}}, \\
a_{N+2}\geq \frac{q_{N+1}a_{N+1}}{q_{N+2}-c_{N+2}}\geq \frac{a_{N}q_{N}q_{N+1}}{(q_{N+2}-c_{N+2})(q_{N+1}-c_{N+1})},
\end{gather*}
and in general,
\begin{equation}\label{aux}
a_{N+k+1} \geq \frac{a_{N}q_{N}q_{N+1} \dots q_{N+k}}{(q_{N+1}-c_{N+1}) \dots (q_{N+k+1}-c_{N+k+1})}, \quad k\geq 0.
\end{equation}
From~\eqref{aux1} and~\eqref{aux} we get
\begin{equation}\label{auxDiverg}
a_{N+k+1}> \frac{a_{N}q_N}{q_{N+k+1}}, \quad k\geq 0.
\end{equation}
Thus
\[
\sum_{k = 0}^{\infty}a_{N+k+1}> a_{N}q_N\sum_{k = 0}^{\infty}\frac{1}{q_{N+k+1}}
\]
and therefore \(\sum a_{n}\) diverges.
From~\eqref{aux}
\[
(q_{N+k+1}-c_{N+k+1})a_{N+k+1}\geq \frac{a_{N}q_{N}q_{N+1} \dots q_{N+k}}{(q_{N+1}-c_{N+1}) \dots (q_{N+k}-c_{N+k})}, \quad k\geq 0,
\]
and applying once again~\eqref{aux1} we obtain that
\[
q_{N+k+1}a_{N+k+1}>(q_{N+k+1}-c_{N+k+1})a_{N+k+1}\geq a_{N}q_{N}>0, \quad k\geq 0.
\]
This last set of inequalities implies that
\[
\lim_{n\to\infty} q_{N+k+1}a_{N+k+1}\neq 0
\quad \textup{ and } \quad
\lim_{k\to\infty} (q_{N+k+1}-c_{N+k+1})a_{N+k+1}\neq 0,
\]
so both series \(\sum q_{n}a_{n}\) and \(\sum (q_{n}-c_{n})a_{n}\) diverge.

Note that from~\eqref{auxDiverg} we obtain that
\begin{equation}\label{auxDiverg1}
c_{N+k+1}a_{N+k+1}> a_{N}q_N\frac{c_{N+k+1}}{q_{N+k+1}}, \quad k\geq 0.
\end{equation}
Therefore, if \(\sum \frac{c_{n}}{q_{n}} \) diverges, then it is clear that
\(\sum c_{n}a_{n}\) diverges.

In order to prove \((ii)\) define
\[
q_{n} = \frac{\sum_{i = 1}^{n}c_{i}a_{i}}{a_{n}}, \quad n\geq 1.
\]
 Clearly, this is a~sequence of positive real numbers that satisfies
\[
q_{n}\frac{a_{n}}{a_{n+1}}-q_{n+1}\leq -c_{n+1}, \quad n\geq 1.
\]
Let us show that \(\sum \frac{1}{q_{n}}\) diverges by concluding that the sequence \(\{s_{k}\}\), defined as \(s_{k} = \sum_{i = 1}^{k}\frac{1}{q_{i}}\), for each \(k\geq 1\), is not a~Cauchy sequence.
Since \(\sum c_{n}a_{n}\) is divergent, given \(m\in \mathbb{N}\) there exists \(k > m\), \(k\in\mathbb{N}\), such that
\begin{equation}\label{aux2}
c_{m}a_{m}+\dots+c_{k}a_{k}>c_{1}a_{1}+\dots+c_{m-1}a_{m-1}.
\end{equation}
Also, from the hypothesis, there exists \(r\geq m\) such that
\begin{equation}\label{aux3}
a_{m}+\dots+a_{r}\geq c_{m}a_{m}+\dots+c_{r}a_{r}.
\end{equation}
Next, we split the proof in two cases: \(k\leq r\) and \(k>r\).

If \(k\leq r\), from~\eqref{aux2} we see that
\begin{equation}\label{aux112}
c_{m}a_{m}+\dots+c_{k}a_{k}+\dots+c_{r}a_{r}\geq c_{m}a_{m}+\dots+c_{k}a_{k}>c_{1}a_{1}+\dots+c_{m-1}a_{m-1}.
\end{equation}
Thus, by~\eqref{aux112} and~\eqref{aux3}
\begin{align*}
\sum_{n = m}^{r}\frac{1}{q_{n}}
{ }&{ }= \frac{a_{m}}{c_{1}a_{1}+ \dots +c_{m}a_{m}}+\dots+\frac{a_{r}}{c_{1}a_{1}+ \dots +c_{r}a_{r}}\\
{ }&{ }\geq \frac{a_{m}+ \dots +a_{r}}{c_{1}a_{1}+ \dots +c_{r}a_{r}}\\
{ }&{ }\geq \frac{c_{m}a_{m}+ \dots +c_{r}a_{r}}{c_{1}a_{1}+ \dots +c_{r}a_{r}}\\
{ }&{ }= \frac{1}{\frac{c_{1}a_{1}+ \dots +c_{m-1}a_{m-1}}{c_{m}a_{m}+ \dots +c_{r}a_{r}}+ 1}\\
{ }&{ }>\frac{1}{2}
\end{align*}
and \(\{s_k\}\) is not a~Cauchy sequence. On the other hand, if \(k>r\) we can use hypothesis again (now applied to \(m_1 = r+1\)) and to obtain \(r_{1}\geq r+1\) such that
\[
a_{r+1}+\dots+a_{r_{1}}\geq c_{r+1}a_{r+1}+\dots+c_{r_{1}}a_{r_1}.
\]
Again, we can use the same argument to conclude that there exists \(r_{2}\geq r_{1}+1\) such that
\[
a_{r_{1}+1}+\dots+a_{r_{2}}\geq c_{r_{1}+1}a_{r_{1}+1}+\dots+c_{r_{2}}a_{r_{2}}.
\]
This procedure can be applied a~finite number of times in order to obtain \(r_{j} \geq k\) for which
\[
a_{r_{(j-1)}+1}+\dots+a_{r_{j}}\geq c_{r_{(j-1)}+1}a_{r_{(j-1)}+1}+\dots+c_{r_{j}}a_{r_{j}}.
\]
Summing up~\eqref{aux3} with all these previous inequalities we obtain that
\[
a_{m}+\dots+a_{r_{j}}\geq c_{m}a_{m}+\dots+c_{r_{j}}a_{r_{j}}
\]
with \(k\leq r_j\).
This reduces the proof to the previous case which we have already proved.
\end{proof}

\section{Some examples and consequences}\label{EC}

The main goal in this section is to present some of the implications of the main results of this paper.

The next three theorems are extensions of the Raabe, Bertrand and Gauss test derived from
Theorem~\ref{thm1} and Theorem~\ref{thm2}. For more information about these tests we refer to~\cite{Ludmila}, \cite{Knopp} and references therein.

Consider the sequences
\[
R_n^{-} = n \frac{a_n}{a_{n+1}}-(n+1)- c_{n+1} \quad \textrm{and} \quad \,R_n^{+} = n \frac{a_n}{a_{n+1}}-(n+1)+ c_{n+1},
\]
for all positive integer \(n\).

\begin{theorem}[Raabe's test]
Let \(\sum c_n a_n\) be a~series of positive terms and suppose
that \(\liminf R_n^{-} = R_1\) and \(\limsup R_n^{-} = R_2\). If
\begin{enumerate}
\item[(i)] \(R_1> 0\), then \(\sum c_n a_n\) converges;
\item[(ii)] \(R_2< 0\) and \(\sum c_n/n \) diverges, then \(\sum c_n a_n\) diverges.
\end{enumerate}
\end{theorem}

\begin{proof}
\begin{enumerate}
\item[(i)] If \(R_1> 0\), then for all \(n\) sufficiently large we have that
\[
n\frac{a_n}{a_{n+1}}-(n+1)- c_{n+1}\geq 0,
\]
hence Theorem~\ref{thm1}, with \(q_n = n\) for all \(n\geq 1\), implies that the series \(\sum c_n a_n\) converges.

\item[(ii)] If \(R_2< 0\), then for all \(n\) sufficiently large
\[
n\frac{a_n}{a_{n+1}}-(n+1)+ c_{n+1}\leq 0.
\]
Again, we have \(q_{n} = n\) for all \(n\geq1\). So, due to the divergence of \(\sum c_n/n \), Theorem~\ref{thm2} implies that \(\sum c_n a_n\) diverges.
\qedhere
\end{enumerate}
\end{proof}

\begin{theorem}[Bertrand's test]
Let \(\sum c_n a_n\) be a~series of positive terms.
\begin{enumerate}
\item[(i)] If
\[
\frac{a_n}{a_{n+1}}> 1+\frac{1}{n}+\frac{\theta_n +c_{n+1}}{n\ln(n)},
\]
for some sequence \(\{\theta_n\}\), such that \(\theta_{n}\geq \theta>1\), for all \(n\geq 1\), then \(\sum c_n a_n\) converges.

\item[(ii)] If
\[
\frac{a_n}{a_{n+1}}\leq 1+\frac{1}{n}+\frac{\theta_n -c_{n+1}}{n\ln(n)},
\]
for some sequence \(\{\theta_n\}\), such that \(\theta_{n}\leq \theta<1\), for all \(n\geq 1\), and \(\sum \frac{c_n}{n\ln(n)} \) diverges, then \(\sum c_n a_n\) diverges.
\end{enumerate}
\end{theorem}

\begin{proof}
\begin{enumerate}
\item[(i)] From the assumption
we get
\[
n\ln(n)\frac{a_n}{a_{n+1}}\geq n\ln(n)+\ln(n)+ c_{n+1}+\theta_n,
\]
for all \(n\) sufficiently large.
That is,
\[
n\ln(n)\frac{a_n}{a_{n+1}}-(n+1)\ln(n+1)\geq (n+1)\ln\left(\frac{n}{n+1}\right)+\theta_n+ c_{n+1},
\]
for all \(n\) sufficiently large.
It follows from the assumption on \(\{\theta_n\}\) that
\[
(n+1)\ln\left(\frac{n}{n+1}\right)+\theta_n> 0,
\]
for all \(n>1\) sufficiently large
 hence
we conclude that
\[
n\ln(n)\frac{a_n}{a_{n+1}}-(n+1)\ln(n+1)> c_{n+1},
\]
for all \(n\) sufficiently large.
Therefore, the convergence of \(\sum c_{n}a_{n}\) follows from an application of Theorem~\ref{thm1}.

\item[(ii)] It suffices to note that
\[
n\ln(n)\frac{a_n}{a_{n+1}}-(n+1)\ln(n+1)\leq (n+1)\ln\left(\frac{n}{n+1}\right)+\theta_n- c_{n+1},
\]
for all \(n\) sufficiently large.
Since \((n+1)\ln\left(\frac{n}{n+1}\right)+\theta_n< 0\) for all \(n>1\) sufficiently large
we obtain
\[
n\ln(n)\frac{a_n}{a_{n+1}}-(n+1)\ln(n+1)< - c_{n+1},
\]
for all \(n\) sufficiently large. The conclusion follows from Theorem~\ref{thm2}.
\qedhere
\end{enumerate}
\end{proof}

\begin{theorem}[Gauss's test]\label{gauss}
Let \(\sum c_n a_n\) be a~series of positive terms, \(\gamma\geq 1\) and \(\{\theta_n\}\) a~bounded sequence of real numbers.
\begin{enumerate}
\item[(i)] Suppose that there exists a~\(\mu \in\mathbb{R}\) such that \(\theta_{n}\geq (1-\mu)n^{\gamma-1}\) holds for all \(n\) sufficiently large. If
\[
\frac{a_n}{a_{n+1}}\geq 1+ \frac{c_{n+1}}{n}+\frac{\mu}{n}+\frac{\theta_n}{n^{\gamma}},
\]
holds for all \(n\) sufficiently large, then
\(\sum c_n a_n\) converges.

\item[(ii)] Suppose that there exists a~\(\mu \in\mathbb{R}\) such that \(\theta_{n}\leq (1-\mu)n^{\gamma-1}\) holds for all \(n\) sufficiently large. If \(\sum c_n/n \) diverges and
\[
\frac{a_n}{a_{n+1}}\leq 1- \frac{c_{n+1}}{n}+\frac{\mu}{n}+\frac{\theta_n}{n^{\gamma}},
\]
for all \(n\) sufficiently large, then
\(\sum c_n a_n\) diverges.
\end{enumerate}
\end{theorem}

\begin{proof}
\begin{enumerate}
\item[(i)] From the assumption we obtain that
\[
n\frac{a_n}{a_{n+1}}-(n+1)\geq c_{n+1}+(\mu-1)+\frac{\theta_n}{n^{\gamma-1}},
\]
 for all \(n\) sufficiently large.
 Taking \(N>0\) such that \(\mu-1+\frac{\theta_n}{n^{\gamma-1}}\geq 0\),
for all \(n>N\), we concude that
\[
n\frac{a_n}{a_{n+1}}-(n+1)\geq c_{n+1},
\]
for all \(n>N\). Therefore, by Theorem~\ref{thm1}, the series \(\sum c_n a_n\) converges.

\item[(ii)] Due to the assumptions on \((ii)\), we have that \(\mu-1 +\frac{\theta_n}{n^{\gamma-1}}<0\) and
\[
n\frac{a_n}{a_{n+1}}-(n+1)\leq -c_{n+1}+(\mu-1)+\frac{\theta_n}{n^{\gamma-1}}\leq -c_{n+1},
\]
 for all \(n\) sufficiently large. The conclusion follows from an application of Theorem~\ref{thm2}.
\qedhere
\end{enumerate}
\end{proof}

Theorem~\ref{thm1} also allows us to provide a~different approach for the well-know Cauchy's condensation test, which we present in the next lemma.

\begin{lemma}\label{CC}\cite[p. 120]{Knopp}(Cauchy's condensation test)
Let \(\{a_n\}\) be a~decreasing sequence of positive numbers.
Then \(\sum a_{n}\) converges if, and only if, \(\sum 2^{n}a_{2^{n}}\) converges.
\end{lemma}

For a~decreasing sequence \(\{a_{n}\}\) of positive real numbers, combining Lemma~\ref{CC} with Theorem~\ref{thm1}, we obtain a~the following characterization of convergence.

\begin{theorem}
Let
\(\sum a_{n}\) be a~series with \(\{a_n\}\) being a~decreasing sequence. Then \(\sum a_{n}\) converges if, and only if,
there exists a~sequence \(\{q_{n}\}\) of positive numbers such that
\[
q_{n}-2 q_{n+1}\geq 2 a_{2^{n+1}},
\]
for all \(n\) sufficiently large.
\end{theorem}
\begin{proof}
By Lemma~\ref{CC}, \(\sum a_n\) coverges if, and only if,
\(\sum 2^n a_{2^n}\) converges. On the other hand, an application of Theorem~\ref{thm1} with \(a_{n} = 2^n\) and \(c_n = a_{2^{n}}\)
show us that \(\sum 2^n a_{2^n}\) converges if, and only if,
there exists a~sequence \(\{q_n\}\) of positive real numbers
such that
\[
q_n-2 q_{n+1}\geq 2 a_{2^{n+1}},
\]
for all \(n\) sufficiently large.
The proof is concluded.
\end{proof}

To close this section of applications we present a~result related to the Olivier's Theorem, which is stated below.
\begin{lemma}[{\cite[p. 124]{Knopp} or \cite{Const}}] \label{Olivier}
Let \(\{a_{n}\}\) be summable decreasing 	sequence of positive real numbers. Then \(\lim n\,a_{n} = 0\).
\end{lemma}

We are going to show that it possible to recover the same Olivier's asymptotic behavior for \(\{a_n\}\) without the decreasigness assumption on \(\{a_n\}\).
Instead of using the monotonicity, we consider an additional assumption on the sequence \(\{q_{n}\}\) (that auxiliary sequence of Theorem~\ref{thm1}). 

\begin{theorem}
Suppose that \(\{a_{n}\}\) is a~sequence of positive numbers. We have that \(\sum a_n\) converges
if, and only if, there exists a~sequence \(\{q_n\}\) of positive numbers such that
\[
q_{n}\frac{n+1}{n}-q_{n+1}\geq (n+1)a_{n+1},
\]
for all \(n\) sufficiently large.
Moreover, if \(\{q_n\}\) satisfies
\[
\lim q_{n}\frac{n+1}{n}-q_{n+1} = 0,
\]
then \(\lim n a_{n} = 0\).
\end{theorem}
\begin{proof}
It is clear that \(\sum a_{n}\) converges if and only if \(\sum \frac{1}{n} na_{n}\) also converges.
From Theorem~\ref{thm1}, with \(a_n = 1/n\) and \(c_n = n a_n\), we can conclude that \(\sum a_n\) converges if, and only if, there exists
a sequence \(\{q_{n}\}\) such that
\[
q_{n}\frac{n+1}{n}-q_{n+1}\geq (n+1)a_{n+1},
\]
for all \(n\) sufficiently large.
Hence, \(\lim n a_n = 0\) certainly occurs when
the sequence \(\{q_n\}\) above is such that
\[
\lim q_{n}\frac{n+1}{n}-q_{n+1} = 0.
\qedhere
\]
\end{proof}

For more information on this asymptotic behavior of summable sequences of positive numbers we refer to~\cite{Lifly}, \cite{Const}, \cite{Salat} and references therein.

\EditInfo{%
	26 February 2019}{%
	17 June 2020}{%
	Karl Dilcher}

\end{paper}